\documentclass{article}
\usepackage{lineno,hyperref}
\modulolinenumbers[5]

\usepackage{amsthm,amsmath,amsfonts,amssymb,graphicx,color,subfigure}

\newtheorem {theorem} {Theorem}

\newtheorem {remark} [theorem]{Remark}

\def\r{{\mathbf{r}}}

\begin{document}


\title{Stability of the regular $n$-gon rotating equilibria with logarithm  interaction}

\author{Anna-Monika Musca\c s\footnote{The Center for Doctoral University Studies, 
University of Oradea,  Oradea, Romania, 
Email: monicaszilagyi@gmail.com} \,\,\,
Daniel Pa\c sca\footnote{Department of Mathematics and Informatics, University of Oradea,  Oradea, Romania, 
Email: dpasca@uoradea.ro} \,\,\,
Cristina Stoica\footnote{Department of Mathematics, 
Wilfrid Laurier University, Waterloo, N2L 3C5, Canada, 
Email: cstoica@wlu.ca}}


\maketitle


\noindent
\textbf{Abstract} 
We study the linear stability of regular $n$-gon rotating equilibria in the $n$-body problem with logarithm interaction. In the presence of a central mass $M$,  linear stability is insured if $M$ is  bounded below and above by constants depending on the number and mass of the (equal) outer $n$ bodies. Moreover, we provide explicit equations of these bounds.
In the absence of a central mass we find that the regular $n$-gon  is linearly stable  for $n =2,3,\ldots 6$ only.

\bigskip

\bigskip



\section{Introduction}

 The logarithm function is the solution of the Laplace equation in two dimensions.  Analogously to the Newtonian potential, one may accept that in a two-dimensional  Euclidean universe the gravitational potential between two mass points $m_1$ and $m_2$ is given by
$
U(\r)=Gm_1m_2 \ln r\,
$
where $\r$ is the vector between the two points and $G$ is a constant.  Despite their  lack of  direct physical relevance, logarithm  potentials  are used in astrophysics  in order to construct models of galaxies that are self-consistent; see   for example, \cite{BiTr87, MESc89}, and more recently,  \cite{BBP07, VWD12}.  
From the standpoint of celestial mechanics,  there are few studies of the  logarithm $n$-body problem. Previous work on the problem may be found in \cite{CaTe11}, where 
    the authors prove  that in the logarithm central force problem collisional solutions may be replaced by transmission trajectories.  The regularisation of the anisotropic case  is investigated in detail  in \cite{StFo03}. 
    In    \cite{Vi07} the authors  prove the existence of periodic solutions. In \cite{MPS24}  the authors study equilibria and stability in the restricted three-body problem. Finally, in a recent paper \cite{SaSt23}, the logarithm central force problem is shown to be  regularisable in the Conley and Easton sense (or \textit{block} regularisable).


 An attractive law,  the logarithm   ``pull"   is weaker at close range than in any law of the form $-1/r^{\alpha},$ $\alpha>0$, but stronger at long range. Since all  trajectories  in the logarithm central force problem  are bounded (and not necessarily periodic),  rotating equilibria (RE) or any other dynamical structures (periodic orbits, invariant tori) perhaps  are stable for a large set of parameters (but this is to be proven). A particular case is that of the regular $n$-gon RE formed by $n$ equal mass-points with or without around a centrally-located mass. As known,  for Newtonian interactions, the regular $n$-gon RE is unstable for $n=2,3\ldots 7$, whereas for $n\geq 7$ linear stability is insured provided there is a  sufficiently heavy central mass \cite{Ro00}.
  In the case of the logarithm interaction we also find that   for large $n$ linear stability is insured provided the 
  the central mass is sufficiently heavy, but moreover, an upper bound is required. (In simple terms,  the centrals mass must be sufficiently heavy to keep the outer masses from escaping but not too heavy, so that falling central mass is prevented.) Specifically, denoting $\mu=m/M$, where $m$ and $M$  the masses of the outer and central bodies, respectively, we obtain that the regular $n$-gon RE is
\begin{itemize}
 \item unstable for  $n=2$;
  \item linearly stable  for  $n=3$ iff $\mu=1$ (i.e. all masses are equal);
  \item linearly stable for  $n=4, 5, 6, 7, 8, 9$  iff 
  $\mu \in [4/(n-1)^2, 1)$;\
  \item linearly stable for $n \geq 10$ even  iff  
  $\mu \in [4/(n-1)^2, 16/(n^2-8n+8)]$;
  \item linearly stable for  $n \geq 11$ odd iff   $\mu \in [4/(n-1)^2, 16/n^2-8n+7)].$

%
\end{itemize} 
In the absence of a central mass, the $n$-gon RE is stable for $n=2,3,\ldots 6$.

For our calculations we adopt the straight-forward method used in \cite{VK07}. The software Matematica proved to be of great help  as it evaluated most of the difficult sums appearing all along.  The paper is organised as follows: in Section 2 
we set up the problem and  find the RE.  In the next section we calculate the linearisation matrix that,  given the symmetry of the regular $n$-gon RE,  displays circulant structures.  Further,  the Hamiltonian nature of the system leads to  a factorisation of the characteristic polynomial   in polynomials of the form $(\lambda^4+a\lambda^2+b).$ In Section 4 we perform all necessary calculation and conclude with  the main Theorem \ref{th: main}, whereas Section 5 concerns the no-central mass case.

\section{Set-up}

Consider the  planar $(n+1)$-body  problem, $n\geq 2,$ with of one large  central body having mass M and $n$ mass points each of mass $m$ orbiting the large body in circular orbits uniformly spaced in a ring of radius $r$.The mutual interaction between any two bodies is given by the logarithm potential. Indices $0$ to $(n-1)$ are  denote the ring masses whereas  index $n$ is used for central body. The coordinates $(x,y)$ of a body are given in complex notation, that is  $z=x+iy$. Thus the equations of motion are:
\begin{equation}\label{e1}
\ddot z_j = - GM \frac{z_j-z_n}{|z_j-z_n|^2} + \sum \limits_{k\neq j,n}Gm\frac{z_k-z_j}{|z_k-z_j|^2}
\end{equation}
Choosing  characteristic scales $(t_0, r_0)$  so that $r_0^2= GM t_0^2$  and denoting $\mu:= m/M$, the equations of motion become
\begin{equation}\label{e2}
\ddot z_j = - \frac{z_j-z_n}{|z_j-z_n|^2} + \sum \limits_{k\neq j,n}\mu\frac{z_k-z_j}{|z_k-z_j|^2}\,.
\end{equation}
A regular $n$-gon relative equilibrium (RE) is a solution of the form 
\begin{align}\label{e3}
&z_j = r e^{i(\omega t + \theta_j)}\,, \,\,\, \text{with}\,\,\,\theta_j:=\frac{2 \pi j}{n} \,,\,\,\, j=0,1,\ldots, (n-1), \,\,\,\\
&z_n = 0\,, \nonumber
\end{align}
for some $\omega \in \mathbb{R}$\,. 
The radius of the $n$-gon  circumcircle is $r>0$ and $ \omega$ is the (uniform) angular velocities of the outer bodies. Note that since the force is attractive there are no equilibrium solutions  (i.e. with $\omega =0$). 
To determine $\omega$, given the symmetry of the problem and that the central body is fixed at the origin,   it is suffices to consider the position of one outer body in interaction with all other. 
From equation \eqref{e3} we have that
%
\begin{equation}\label{e4}
\ddot z_j = - \omega^2 z_j\,\,\,\,\text{for}\,\,\,j=0,1,\ldots, (n-1).
\end{equation}
Without loosing generality, for $j=0$ from (\ref{e3}) we have that
\begin{equation}\label{e5}
z_k - z_0 = r  e^{i(\omega t)} e^{i\left(\frac{\theta_k}{2}\right)}  2 i \sin\left(\frac{\theta_k}{2}\right),
\end{equation}
and so
\begin{equation}\label{e6}
|z_k - z_0| = 2 r \sin\left(\frac{\theta_k}{2}\right).
\end{equation}
Substituting (\ref{e5}) and (\ref{e6}) into (\ref{e2}) , and using \eqref{e4}, we obtain:
\begin{equation}\label{e7}
\omega^2 =  \frac{1}{r^2} + \frac{\mu}{2 r^2} \sum_{k=1}^{n-1} \Big(1 - \frac{ i \cos \frac{\theta_k}{2}}{\sin \frac{\theta_k}{2}}\Big).
\end{equation}
Elementary calculations show that the imaginary  part of the sum above is zero and so we obtain the relation between the angular velocity $\omega$ and the radius of a regular polygon solution:
\begin{equation}\label{e8_n}
\omega^2 = \frac{1}{r^2} + \frac{\mu(n-1)}{2r^2}.
\end{equation}
or
\begin{equation}\label{e8}
r^2 \omega^2=1+\frac{\mu(n-1)}{2}\,.
\end{equation}

\begin{remark}
In the Newtonian case, we have
\begin{equation}\label{eq-RE-Newton}
r^3 \omega^2= 1+ \frac{\mu}{2}\sum 
\limits_{k=1}^{n-1} \frac{1}{2\sin (\theta_k/2)}
\end{equation}

\end{remark}

\section{Linearization}

%
%
We start by applying the change of variables $z_j \to w_j$  given by
\begin{equation}\label{st2}
w_j = u_j + i v_j = e^{-i(\omega t + 2\pi j/n)}z_j\,,\,\,\,\,j=0,1,\ldots,n-1,n\,.
\end{equation}
so in the new coordinates the configuration of the relative equilibrium  reads
\begin{equation}\label{st7}
w_{j,e} =
\begin{cases}
r; \,\,\, \text{for}\,\,\, k=0, 1,\ldots, (n-1)\\
0; \,\,\, \text{for}\,\,\,  k=n.
\end{cases}
\end{equation}
%
%
Differentiating (\ref{st2}) twice, we get
\begin{equation}\label{st3}
\ddot w_{j} = \omega^2 w_{j}-2i\omega \dot w_{j}+ e^{-i(\omega t+ \theta_{j})}\ddot z_{j}
\end{equation}
from where, using (\ref{e1}) we have
\begin{equation}\label{st4}
\ddot w_{j} = \omega^{2}w_{j}-2i\omega \dot w_{j} + \sum_{k \neq j} m_{k} \frac{\varepsilon_{k,j}}{|\varepsilon_{k,j}|^2}, \,\,\,   j=0, 1,\ldots, (n-1)
\end{equation}
where
\begin{equation}\label{st5}
\varepsilon_{k,j} := w_{k}e^{i\theta_{k-j}}-w_{j},
\end{equation}
and
\begin{equation}\label{st6}
m_{k}:=
\begin{cases}
\mu; \,\,\, \text{for} \,\,\, k=0, 1,\ldots, (n-1)\\
1; \,\,\, \text{for} \,\,\, k=n.
\end{cases}
\end{equation}
Calculating the variations $\delta w_j(t)$  about  $w_{j,e}$, $j=0,1,\ldots, n$ we obtain
%
%
\begin{equation}\label{st8}
\delta\ddot w_{j}=\omega^{2}\delta w_{j}-2i\omega \delta\dot w_{j} + \sum_{k \neq j} m_{k} \delta \Big(\frac{\varepsilon_{k,j}}{|\varepsilon_{k,j}|^2}\Big)\,.
\end{equation}
One may verify that
\begin{equation}\label{st9}
\delta\Big(\frac{\varepsilon_{k,j}}{|\varepsilon_{k,j}|^2}\Big) = - \frac{\varepsilon_{k,j}^{2}\delta\overline{\varepsilon}_{k,j}}{|\varepsilon_{k,j}|^4}.
\end{equation}
For $k\neq n$ we have
\begin{equation}
\label{st10}
-\frac{\varepsilon_{k,j}^{2}\delta\overline\varepsilon_{k,j}}{|\varepsilon_{k,j}|^4} = \frac{\delta\overline w_{k}-e^{-i\theta_{k-j}}\delta\overline w_{j}}{4r^2\sin^2\frac{\theta_{k-j}}{2}}
\end{equation}
whereas for  for $k=n$ we get
\begin{equation}\label{st11}
-\frac{\varepsilon_{n,j}^{2}\delta\overline\varepsilon_{n,j}}{|\varepsilon_{k,j}|^4}
 =-\frac{1}{r^2}(e^{i\theta_{j}}\delta\overline w_{n}-\delta\overline w_{j}).
\end{equation}
Since the centre of mass is fixed at the origin, 
 conservation of momentum implies that
\begin{equation}\label{st12}
m \sum_{k\neq n} \delta z_k + M \delta z_n =0,
\end{equation}
and hence
\begin{equation}\label{st13}
\delta z_n = - \mu \sum_{k\neq n} \delta z_k.
\end{equation}
Using  the definition (\ref{st2}) of the $w_k$in terms of  $z_k$, it follows that
\begin{equation}\label{st14}
e^{-i \theta_j}\delta w_n = - \mu \sum_{k\neq n} e^{i \theta_{k-j}}\delta w_k.
\end{equation}
Making this substitution for $e^{-i \theta_j}\delta w_n$ and an 
%
\begin{equation}\label{st15}
\begin{split}
\delta \ddot w_{j} =\omega^{2}\delta w_{j}+\frac{\mu}{r^2}\delta\overline w_{j}+ &\frac{1}{r^2}\delta\overline w_{j}-\frac{\mu}{4r^2}\Big(\sum_{k \neq j,n}\frac{e^{i\theta_{k-j}}}{\sin^2\frac{|\theta_{k-j}|}{2}}\Big)
\delta\overline w_{j}\\
& -2i\omega\delta\dot w_{j}+\frac{\mu}{r^2}\sum_{k \neq j,n}e^{-i\theta_{k-j}}\delta\overline w_{k}
+\frac{\mu}{4r^2}\sum_{k \neq j,n}\frac{\delta \overline w_{k}}{\sin^2\frac{|\theta_{k-j}|}{2}},
\end{split}
\end{equation}
and its conjugate
\begin{equation}\label{st16}
\begin{split}
\delta \ddot {\overline w_{j}}=\omega^{2}\delta \overline w_{j}+\frac{\mu}{r^2}\delta w_{j}+ &\frac{1}{r^2}\delta w_{j}-\frac{\mu}{4r^2}\Big(\sum_{k \neq j,n}\frac{e^{-i\theta_{k-j}}}{\sin^2\frac{|\theta_{k-j}|}{2}}\Big)
\delta w_{j}\\
&+2i\omega\delta\dot {\overline w_{j}}+\frac{\mu}{r^2}\sum_{k \neq j,n}e^{i\theta_{k-j}}\delta w_{k}
+\frac{\mu}{4r^2}\sum_{k \neq j,n}\frac{\delta w_{k}}{\sin^2\frac{|\theta_{k-j}|}{2}}\,.
\end{split}
\end{equation}
One may verify  that
\begin{equation}\label{st17}
\begin{split}
\sum_{k \neq j,n}\frac{e^{\pm i\theta_{k-j}}}{\sin^2\frac{|\theta_{k-j}|}{2}}&=\sum_{k \neq j,n}\frac{\cos{\theta_{k-j}} \pm i\sin{\theta_{k-j}}}{\sin^2\frac{|\theta_{k-j}|}{2}}=
\sum_{k \neq j,n}\frac{1-2\sin^2\frac{\theta_{k-j}}{2}\pm 2i\sin\frac{\theta_{k-j}}{2}\cos\frac{\theta_{k-j}}{2}}{\sin^2\frac{|\theta_{k-j}|}{2}}=\\
&=\sum_{k \neq j,n}\frac{1}{\sin^2\frac{|\theta_{k-j}|}{2}}-2(n-1)=\frac{n^2-1}{3}-2(n-1)=\frac{(n-1)(n-5)}{3},
\end{split}
\end{equation}
and
\begin{equation}\label{st18}
\begin{split}
\frac{\mu}{4r^2}\Big(\sum_{k \neq j,n}\frac{e^{\pm i\theta_{k-j}}}{\sin^2\frac{|\theta_{k-j}|}{2}}\Big) &=\frac{\mu}{4}\frac{\omega^2}{1+\frac{\mu(n-1)}{2}}\frac{(n-1)(n-5)}{3}
=\frac{\mu\omega^2(n-1)(n-5)}{6[2+\mu(n-1)]} =: a.
\end{split}
\end{equation}
So, finally we can write the equations (\ref{st15}) and (\ref{st16}) on the following form
\begin{equation}\label{st19}
\delta \ddot w_{j} =\omega^{2}\delta w_{j}+\Big( \frac{\mu +1}{r^2} -a\Big) \delta\overline w_{j} -2i\omega\delta\dot w_{j}+
\frac{\mu}{r^2} \sum_{k \neq j,n} \Big( e^{-i\theta_{k-j}} + \frac{1}{4\sin^2\frac{\theta_{k-j}}{2}}\Big)
\delta\overline w_{j},
\end{equation}
\begin{equation}\label{st20}
\delta \ddot {\overline w}_{j} =\omega^{2}\delta \overline w_{j}+\Big( \frac{\mu +1}{r^2} -a\Big) \delta w_{j} +2i\omega\delta\dot{\overline w}_{j}+
\frac{\mu}{r^2} \sum_{k \neq j,n} \Big( e^{i\theta_{k-j}} + \frac{1}{4\sin^2\frac{\theta_{k-j}}{2}}\Big)
\delta w_{j}.
\end{equation}
%
%
%
%
Let $W_{j}$ denote a shorthand for the vector 
$\displaystyle{
\left[
\begin{array}{c}
w_{j}\ \\
  \overline w_{j}
\end{array}
\right] \in \mathbb{C}^2
}$. 
In this notation, we see that (\ref{st19}) and (\ref{st20}) can be written as
\begin{equation}\label{cm1}
\frac{d}{dt}
\begin{bmatrix}
\delta W_{0}\\
\delta W_{1}\\
\vdots \\
\delta W_{n-1}\\
\delta\dot W_{0}\\
\delta\dot W_{1}\\
\vdots\\
\delta\dot W_{n-1}
\end{bmatrix}
=
\begin{bmatrix}
\begin{array}{c | c}
   0_{2n,2n}
    &
    \begin{array}{c c c c}
     I_{2,2} & 0_{2,2} & \ldots & 0_{2,2}\\
     0_{2,2} & I_{2,2} & \ldots & 0_{2,2}\\
     \ldots & \ldots & \ldots & \ldots\\
     0_{2,2} & 0_{2,2} &\ldots & I_{2,2}
    \end{array} \\
  \hline
  \begin{array}{c c c c}
     D & N_{1} & \ldots & N_{n-1}\\
     N_{1} & D & \ldots & N_{n-2}\\
     \ldots & \ldots & \ldots & \ldots\\
     N_{1} & N_{2} &\ldots & D
  \end{array}&
  \begin{array}{c c c c}
     \Omega & 0_{2,2} & \ldots & 0_{2,2}\\
     0_{2,2} & \Omega & \ldots & 0_{2,2}\\
     \ldots & \ldots & \ldots & \ldots\\
     0_{2,2} & 0_{2,2} &\ldots & \Omega
  \end{array}
 \end{array}
 \end{bmatrix}
\begin{bmatrix}
\delta W_{0}\\
\delta W_{1}\\
\vdots \\
\delta W_{n-1}\\
\delta\dot W_{0}\\
\delta\dot W_{1}\\
\vdots\\
\delta\dot W_{n-1}
\end{bmatrix}
\end{equation}
where $D$, $\Omega$ and the $N_k$'s are $2\times 2$ complex matrices given by
\begin{equation}
D=\omega^2
\begin{bmatrix}\label{cm2}
\begin{array}{c c}
     1 & 0\\
     0 & 1
  \end{array}
\end{bmatrix}
+\Big(\frac{\mu+1}{r^2}-a\Big)
\begin{bmatrix}
\begin{array}{c c}
     0 & 1\\
     1 & 0
  \end{array}
\end{bmatrix}
\end{equation}
\begin{equation}\label{cm3}
N_{k}=\frac{\mu}{r^2}
\begin{bmatrix}
\begin{array}{c c}
     0 & e^{-i\theta_{k}}+\frac{1}{4\sin^2\frac{|\theta_{k}|}{2}}\\
     e^{i\theta_{k}}+\frac{1}{4\sin^2\frac{|\theta_{k}|}{2}} & 0
  \end{array}
\end{bmatrix}
\end{equation}
\begin{equation}\label{cm4}
\Omega = 2i \omega
\begin{bmatrix}
\begin{array}{c c}
-1 & 0\\
0 & 1
\end{array}
\end{bmatrix}
\end{equation}
with $a$ defined in \eqref{st18}. 
 The (complex) eigenvalues $\lambda$  of the linearization  matrix in the right hand side of (\ref{cm1})
%
are determined by solving
\begin{equation}\label{cm5}
\begin{bmatrix}
\begin{array}{c | c}
   0_{2n,2n}
    &
    \begin{array}{c c c c}
     I_{2,2} & 0_{2,2} & \ldots & 0_{2,2}\\
     0_{2,2} & I_{2,2} & \ldots & 0_{2,2}\\
     \ldots & \ldots & \ldots & \ldots\\
     0_{2,2} & 0_{2,2} &\ldots & I_{2,2}
    \end{array} \\
  \hline
  \begin{array}{c c c c}
     D & N_{1} & \ldots & N_{n-1}\\
     N_{n-1} & D & \ldots & N_{n-2}\\
     \ldots & \ldots & \ldots & \ldots\\
     N_{1} & N_{2} &\ldots & D
  \end{array}&
  \begin{array}{c c c c}
     \Omega & 0_{2,2} & \ldots & 0_{2,2}\\
     0_{2,2} & \Omega & \ldots & 0_{2,2}\\
     \ldots & \ldots & \ldots & \ldots\\
     0_{2,2} & 0_{2,2} &\ldots & \Omega
  \end{array}
 \end{array}
 \end{bmatrix}
\begin{bmatrix}
\delta W_{0}\\
\delta W_{1}\\
\vdots \\
\delta W_{n-1}\\
\delta\dot W_{0}\\
\delta\dot W_{1}\\
\vdots\\
\delta\dot W_{n-1}
\end{bmatrix}
=\lambda
\begin{bmatrix}
\delta W_{0}\\
\delta W_{1}\\
\vdots \\
\delta W_{n-1}\\
\delta\dot W_{0}\\
\delta\dot W_{1}\\
\vdots\\
\delta\dot W_{n-1}
\end{bmatrix}.
\end{equation}
Using that
\begin{equation}\label{cm6}
\begin{bmatrix}
\delta\dot W_{0}\\
\delta\dot W_{1}\\
\vdots\\
\delta\dot W_{n-1}
\end{bmatrix}
=\lambda
\begin{bmatrix}
\delta W_{0}\\
\delta W_{1}\\
\vdots \\
\delta W_{n-1}\\
\end{bmatrix}
\end{equation}
we have
\begin{equation}\label{cm7}
\begin{bmatrix}
  \begin{array}{c c c c}
     D & N_{1} & \ldots & N_{n-1}\\
     N_{n-1} & D & \ldots & N_{n-2}\\
     \ldots & \ldots & \ldots & \ldots\\
     N_{1} & N_{2} &\ldots & D
  \end{array}
\end{bmatrix}
\begin{bmatrix}
\delta W_{0}\\
\delta W_{1}\\
\vdots \\
\delta W_{n-1}\\
\end{bmatrix}
+\lambda
\begin{bmatrix}
\begin{array}{c c c c}
     \Omega & 0_{2,2} & \ldots & 0_{2,2}\\
     0_{2,2} & \Omega & \ldots & 0_{2,2}\\
     \ldots & \ldots & \ldots & \ldots\\
     0_{2,2} & 0_{2,2} &\ldots & \Omega
  \end{array}
\end{bmatrix}
\begin{bmatrix}
\delta W_{0}\\
\delta W_{1}\\
\vdots \\
\delta W_{n-1}\\
\end{bmatrix}
=
\lambda^2
\begin{bmatrix}
\delta W_{0}\\
\delta W_{1}\\
\vdots \\
\delta W_{n-1}\\
\end{bmatrix}
\end{equation}
Further, for  the block circulant matrix \cite{D94} we implement the ansatz
%
%
\begin{equation}\label{cm8}
\begin{bmatrix}
\delta W_{0}\\
\delta W_{1}\\
\vdots \\
\delta W_{n-1}\\
\end{bmatrix}
=
\begin{bmatrix}
\xi\\
\rho\xi\\
\vdots \\
\rho^{n-1}\xi\\
\end{bmatrix}
.
\end{equation}
where  $ \rho $ denote an $ n $-th root of unity (i.e., $ \rho = e^{2\pi ij/n} $ for some $ j = 0, 1, \ldots, n-1$) and  $\xi$ an arbitrary complex 2-vector.  
Substituting this into (\ref{cm7}), we obtain
%
\begin{equation}\label{cm9}
(D + \rho N_{1}+ \ldots + \rho^{n-1}N_{n-1})\xi + \lambda \Omega \xi = \lambda^2 \xi
\end{equation}
and so we have to solve
\begin{equation}\label{cm10}
\det(D + \rho N_{1}+ \ldots + \rho^{n-1}N_{n-1} + \lambda \Omega - \lambda^2 I) = 0.
\end{equation}
(Note that to each root of unity corresponds four ``$\lambda$" roots.)
%
%
%
Using equation (\ref{cm3}) we have
\begin{equation}\nonumber
\sum_{k=1}^{n-1} \rho_j^k N_{k}=\frac{\mu}{r^2}
\begin{bmatrix}
\begin{array}{c c}
     0 & \sum_{k=1}^{n-1} \Big[ e^{i(j-1)\theta_{k}}+\frac{\rho_j^k}{4\sin^2\frac{|\theta_{k}|}{2}}\Big]\\
     \sum_{k=1}^{n-1} \Big[e^{i(j+1)\theta_{k}}+\frac{\rho_j^k}{4\sin^2\frac{|\theta_{k}|}{2}}\Big] & 0
  \end{array}
\end{bmatrix} =
\end{equation}
\begin{equation}\label{ee1}
=\frac{\mu}{r^2}
\begin{bmatrix}
\begin{array}{c c}
     0 & -1 + n \delta_{j=1} + C_j\\
     -1 + n \delta_{j=n-1} + C_j & 0
  \end{array}
\end{bmatrix},
\end{equation}
where $\delta_{j=k}$ denotes the Kronecker delta  and
\begin{equation}\label{ee3}
C_j = \frac{1}{4}\sum_{k=1}^{n-1} \frac{\rho_j^k}{\sin^2\frac{|\theta_{k}|}{2}}.
\end{equation}
%
%
Thus
\begin{equation}\nonumber
\det \Big(D + \sum_{k=1}^{n-1} \rho^k N_{k} + \lambda \Omega - \lambda^2 I \Big) =
\begin{vmatrix}
\begin{array}{c c}
\omega^2 - 2 i \omega \lambda - \lambda^2 & \frac{1}{r^2} - a + \frac{\mu}{r^2} (n \delta_{j=1} + C_j ) \\ \\
  \frac{1}{r^2} - a + \frac{\mu}{r^2} (n \delta_{j=1} + C_j )    & \omega^2 + 2 i \omega \lambda - \lambda^2
  \end{array}
\end{vmatrix} =
\end{equation}
\begin{equation}\label{s1}
= \lambda^4 + 2 \omega^2 \lambda^2 + \omega^4 - \Big[ \frac{1}{r^2} - a + \frac{\mu}{r^2} \Big(n \delta_{j=1} + C_j\Big) \Big]
 \Big[ \frac{1}{r^2} - a + \frac{\mu}{r^2} \Big( n \delta_{j=1} + C_j\Big) \Big]\,.
\end{equation}
%
In conclusion, the eigenvalues  are  roots of the quartic equations
\begin{equation}
 \lambda^4 + 2 \omega^2 \lambda^2 + \omega^4 - \Big[ \frac{1}{r^2} - a + \frac{\mu}{r^2} \Big(n \delta_{j=1} + C_j\Big) \Big]
 \Big[ \frac{1}{r^2} - a + \frac{\mu}{r^2} \Big( n \delta_{j=1} + C_j\Big) \Big]
 \label{eq: main_lin}
\end{equation}
for $j=0,1,2, \ldots (n-1).$

\section{Stability for the regular $n$-ring with a central mass}

Given the Hamiltonian nature of the system, the RE is linearly stable if all the eigenvalues are purely imaginary \cite{MHO09}. Denoting $y=\lambda^2$, the equations \eqref{eq: main_lin} written in the form
\[
y^2+A_j y+B_j=0\,
\]
have all roots purely imaginary if and only if both roots $y_1$ and $y_2$ of the above are real and negative. In its turn, provided the roots are real,   entails  that the sum $S_j = y_1 + y_2$ and  the  product $P_j = y_1 y_2$ must be  negative and positive, respectively  (that is $S_j<0$ and $P_j>0)$ for $j=0,1,2, \ldots (n-1).$

\subsubsection{Analysis for $j=0$.}
For $j=0$, since $\rho_0 =1$, we calculate
\begin{equation}\label{s8}
C_0 = \frac{1}{4}\sum_{k=1}^{n-1} \frac{1}{\sin^2 \frac{\pi k}{n}} = \frac{n^2 -1}{12}.
\end{equation}
Thus
\begin{equation}\label{s9}
P_0 = \omega^4 \Big[ 1 - \Big( \frac{12 -\mu (n-1)(n-5) + \mu (n^2 -1)}{6(2+\mu (n-1))}\Big)^2 \Big] = 0,
\end{equation}
and calculating  the eigenvalues we find  
\[\lambda_{1,2}=0\,,\,\,\,\,\,\lambda_{3,4}=\pm i \omega \sqrt{2}\,.\]

\subsection{Analysis for  $j=1$ and $j=n-1$.}

For $j=1$ we have
\begin{equation}\label{s10}
\lambda^4 + 2 \omega^2 \lambda^2 + \omega^4 - \Big( \frac{1}{r^2} - a + \frac{\mu}{r^2} C_1 \Big)  \Big( \frac{1}{r^2} - a +
\frac{\mu}{r^2} C_1 + \frac{\mu}{r^2} n \Big)= 0,
\end{equation}
and
\begin{equation}\label{s11}
\Delta_1 = \Big( \frac{1}{r^2} - a + \frac{\mu}{r^2} C_1 \Big)^2 + \frac{\mu}{r^2} n \Big( \frac{1}{r^2} - a + \frac{\mu}{r^2} C_1 \Big) =
4 \omega^4 \frac{\mu n +1}{(2+\mu (n-1))^2} >0,
\end{equation}
where $\rho_1 = e^{2\pi i /n}$ and
\begin{equation}\label{s12}
C_1 = \frac{1}{4}\sum_{k=1}^{n-1} \frac{e^{\frac{2\pi i k}{n}}}{\sin^2 \frac{\pi k}{n}} = \frac{(n-1)(n-5)}{12}\,.
\end{equation}
Thus
\begin{equation}\label{s13}
P_1 = \omega^4 - \Delta_1 = \omega^4 \mu \frac{\mu (n-1)^2 - 4}{(2+\mu (n-1))^2}.
\end{equation}
$P_1$ is positive iff $\mu \geq \frac{4}{(n-1)^2} $.

\begin{remark}\label{coro1}
For $n=2$ we get $\mu \geq 4$ which is impossible since $\mu \in (0,1)$.
\end{remark}

For $j=n-1$ we have
\begin{equation}\label{s14}
\lambda^4 + 2 \omega^2 \lambda^2 + \omega^4 - \Big( \frac{1}{r^2} - a + \frac{\mu}{r^2} C_{n-1} + \frac{\mu}{r^2} n \Big)  \Big( \frac{1}{r^2} - a +
\frac{\mu}{r^2} C_{n-1} \Big)= 0,
\end{equation}
and
\begin{equation}\label{s15}
\Delta_{n-1} = \Big( \frac{1}{r^2} - a + \frac{\mu}{r^2} C_{n-1} \Big)^2 + \frac{\mu}{r^2} n \Big( \frac{1}{r^2} - a + \frac{\mu}{r^2} C_{n-1} \Big).
\end{equation}
For $j=n-1$ we have $\rho_{n-1} = e^{2\pi i (n-1)/n} = e^{-2\pi i/n}$ and
\begin{equation}\label{s16}
C_{n-1} = \frac{1}{4}\sum_{k=1}^{n-1} \frac{e^{\frac{2\pi i k (n-1)}{n}}}{\sin^2 \frac{\pi k}{n}} = C_1,
\end{equation}
thus $\Delta_{n-1} = \Delta_1 >0$ and $P_{n-1} = P_1$.

\begin{remark}\label{coro2}
For $n=3$ we have $P_1 = P_2$ positive iff $\mu =1$, so the equilateral triangle case is linearly
stable iff the central mass are equal with the masses situated in the triangle vertices.
\end{remark}


\subsection{Analysis for $j\neq 0$,   $j \neq 1$ and $j \neq n-1$.}


In this case the equation \eqref{eq: main_lin}  becomes
\begin{equation}\label{s2}
\lambda^4 + 2 \omega^2 \lambda^2 + \omega^4 - \Big( \frac{1}{r^2} - a + \frac{\mu}{r^2} C_j \Big)^2 = 0.
\end{equation}
Denoting $y=\lambda^2$, the equation above writes
\begin{equation}\label{s3}
y^2 + 2 \omega^2 y + \omega^4 - \Big( \frac{1}{r^2} - a + \frac{\mu}{r^2} C_j \Big)^2 = 0.
\end{equation}
%
%
%
%
The discriminant of the above quadratic is:
\begin{equation}\label{s4}
\Delta_j = \Big( \frac{1}{r^2} - a + \frac{\mu}{r^2} C_j \Big)^2 \geq 0
\end{equation}
and so the equation always has real roots. 
We have
\begin{equation}\label{s5}
S_j = y_1 + y_2 = - 2 \omega^2 \leq 0\,.
\end{equation}
%
For the product $P_j$ we have
 %
\begin{equation}\label{s6}
P_j = y_1  y_2 = \omega^4 - \Big( \frac{1}{r^2} - a + \frac{\mu}{r^2} C_j \Big)^2
\end{equation}
from where, taking into account (\ref{st18}) we get
\begin{equation}\label{s7}
P_j = \omega^4 \Big[ 1 - \Big( \frac{12 -\mu (n-1)(n-5) + 12 \mu C_j}{6(2+\mu (n-1))}\Big)^2 \Big]\,.
\end{equation}
Since
\begin{equation}\label{s54}
C_j = \frac{1}{4}\sum_{k=1}^{n-1} \frac{\rho_j^k}{\sin^2\frac{\theta_{k}}{2}}, \,\,\,\, \rho_j=e^{\frac{2\pi i j}{n}}, \,\,\, j=0,1,\ldots n-1, \,\,\,\, \theta_k = \frac{2\pi k}{n},
\end{equation}
and
\begin{equation}\label{s55}
\frac{\sin \frac{2\pi k j}{n}}{\sin^2\frac{\pi k}{n}} = - \frac{\sin \frac{2\pi (n-k) j}{n}}{\sin^2\frac{\pi (n-k)}{n}}, 
\end{equation}
we obtain that the imaginary part of $C_j$ is zero.
Using Mathematica, we further obtain that
\begin{equation}\label{s56}
C_j = \frac{1}{4}\sum_{k=1}^{n-1} \frac{\cos \frac{2\pi k j}{n}}{\sin^2\frac{ \pi k}{n}} = \frac{1}{12} (n^2 - 6nj + 6j^2 -1), \,\,\,\,j=0,1,\ldots (n-1). 
\end{equation}
It follows that 
\begin{equation}\label{s57}
P_j = \omega^4 \Big[ 1 - \Big( \frac{12 -\mu (n-1)(n-5) + 12 \mu C_j}{6(2+\mu (n-1))}\Big)^2 \Big] = \omega^4 
\frac{\mu j (n-j) \Big[ \mu \Big( j^2 - jn + 2(n-1) \Big) +4 \Big]}{\Big[ \mu (n-1) +2 \Big]^2}.
\end{equation}
So, 
$P_j \geq 0$  when $\mu \Big( j^2 - jn + 2(n-1) \Big) +4 \geq 0.$
We now look at the sign of the parabola $j^2 - jn + 2(n-1)$ for $j=2,3\ldots n-2$. The discriminant $\Delta = n^2 -8n +8$ has the roots $n_{1,2} = 4\pm 2\sqrt{2}$. Thus if $n=2,3,4,5,6$ then $\Delta <0$ and therefore $j^2 - jn + 2(n-1) >0$ for any $j$ and consequently $P_j \geq 0$.  Further
\begin{itemize}
\item for even $n \geq 8$, the minimum of the $j^2 - jn + 2(n-1)$ is for $j=n/2$ and its equal with $-(n^2-8n+8)/4<0$, so $P_j \geq 0$ if $\mu \leq 16/(n^2-8n +8)$. But for $n=8$ this is satisfies because $\mu \leq 1$. For even $n\geq 10$ we will have $P_j \geq 0$ iff $\mu \leq 16/(n^2-8n +8);$

\item for odd $n \geq 9$, the minimum of the $j^2 - jn + 2(n-1)$ is for $j=(n-1)/2$ and for $j=(n+1)/2$ and its equal with $-(n-1)(n-7)/4<0$, so $P_j \geq 0$ if $\mu \leq 16/(n-1)(n-7)$. But for $n=9$ this is satisfies because $\mu \leq 1$. For odd $n\geq 11$ we will have $P_j \geq 0$ iff $\mu \leq 16/(n-1)(n-7)\,.$

\end{itemize}
Recalling the calculations  in the cases $j=0,$ $j=1,$ and $j=n-1$,  we have proven 

\begin{theorem} \label{th: main}Consider the  regular $n$-gon with a central mass relative equilibrium  in the logarithm $n$-body problem and let  $\mu=m/M$ where $m$ and $M$ are the outer and the central masses, respectively. Then the regular $n$-ring is

\begin{itemize}
  \item unstable for  $n=2$ 
  \item linearly stable  for  $n=3$ if $\mu=1$ (i.e. all masses are equal)
  \item linearly stable for  $n=4, 5, 6, 7, 8, 9$  iff $\mu \in [4/(n-1)^2, 1)$.
  \item linearly stable for $n \geq 10$ even  iff $\mu \in [4/(n-1)^2, 16/(n^2-8n+8)]$.
  \item linearly stable for  $n \geq 11$ odd iff $\mu \in [4/(n-1)^2, 16/(n-1)(n-7)]$.
\end{itemize} 
\end{theorem}

\section{Stability for the regular $n$-ring without a central mass}


In the absence of a central mass, relation \eqref{e8} becomes
\begin{equation}\label{wcm1}
r^2\omega^2 =  \frac{n-1}{2}.
\end{equation}
and the variations' equations are
\begin{equation}\label{wcm2}
\delta \ddot w_{j} = \omega^{2}\delta w_{j} - 2i \omega \delta \dot w_{j} +
\sum_{k \neq j,n}  \frac{\delta \overline w_k - e^{i\theta_{k-j}} \delta \overline w_j }{4r^2 \sin^2\frac{\theta_{k-j}}{2}},
\end{equation}
\begin{equation}\label{wcm3}
\delta \ddot {\overline w_{j}} = \omega^{2}\delta \overline w_{j} + 2i \omega \delta \dot {\overline w_{j}} +
\sum_{k \neq j,n}  \frac{\delta w_k - e^{-i\theta_{k-j}} \delta w_j }{4r^2 \sin^2 \frac{\theta_{k-j}}{2}}.
\end{equation}
In the matrix form the system above reads:
\begin{equation}\label{wcm4}
\frac{d}{dt}
\begin{bmatrix}
\delta W_{0}\\
\delta W_{1}\\
\vdots \\
\delta W_{n-1}\\
\delta\dot W_{0}\\
\delta\dot W_{1}\\
\vdots\\
\delta\dot W_{n-1}
\end{bmatrix}
=
\begin{bmatrix}
\begin{array}{c | c}
   0_{2n,2n}
    &
    \begin{array}{c c c c}
     I_{2,2} & 0_{2,2} & \ldots & 0_{2,2}\\
     0_{2,2} & I_{2,2} & \ldots & 0_{2,2}\\
     \ldots & \ldots & \ldots & \ldots\\
     0_{2,2} & 0_{2,2} &\ldots & I_{2,2}
    \end{array} \\
  \hline
  \begin{array}{c c c c}
     D & N_{1} & \ldots & N_{n-1}\\
     N_{1} & D & \ldots & N_{n-2}\\
     \ldots & \ldots & \ldots & \ldots\\
     N_{1} & N_{2} &\ldots & D
  \end{array}&
  \begin{array}{c c c c}
     \Omega & 0_{2,2} & \ldots & 0_{2,2}\\
     0_{2,2} & \Omega & \ldots & 0_{2,2}\\
     \ldots & \ldots & \ldots & \ldots\\
     0_{2,2} & 0_{2,2} &\ldots & \Omega
  \end{array}
 \end{array}
 \end{bmatrix}
\begin{bmatrix}
\delta W_{0}\\
\delta W_{1}\\
\vdots \\
\delta W_{n-1}\\
\delta\dot W_{0}\\
\delta\dot W_{1}\\
\vdots\\
\delta\dot W_{n-1}
\end{bmatrix}
\end{equation}
where $D$, $\Omega$ and the $N_k$'s are the $2\times 2$ matrices:
\begin{equation}
D=\omega^2
\begin{bmatrix}\label{wcm5}
\begin{array}{c c}
     1 & 0\\
     0 & 1
  \end{array}
\,,\,\,\,\,\,
\begin{array}{c c}
     0 & b\\
     b & 0
  \end{array}
\end{bmatrix}
\,,\,\,\,
N_{k} = \frac{1}{4\sin^2\frac{\theta_{k}}{2}}
\begin{bmatrix}
\begin{array}{c c}
     0 & 1 \\
     1 & 0
  \end{array}
\end{bmatrix}
\,,\,\,\,\,
\Omega = 2i \omega
\begin{bmatrix}
\begin{array}{c c}
-1 & 0\\
0 & 1
\end{array}
\end{bmatrix}
\end{equation}
with
$$
b := \sum_{k \neq j,n}  \frac{ e^{-i\theta_{k-j}}}{4r^2 \sin^2\frac{\theta_{k-j}}{2}} = \frac{\omega^2 (n-5)}{6}.
$$
Similar calculations as in the previous section lead to the correspondent  stability equation \eqref{eq: main_lin}:
%
\begin{equation}\label{wcm8}
\lambda^4 + 2 \omega^2 \lambda^2 + \omega^4 - \Big( \frac{1}{r^2} C_j - b \Big)^2 = 0.
\end{equation}
Denoting $y=\lambda^2$, the  above reads
\begin{equation}\label{wcm9}
y^2 + 2 \omega^2 y + \omega^4 - \Big( \frac{1}{r^2} C_j - b \Big)^2 = 0.
\end{equation}
%
%
%
Its discriminant is
\begin{equation}\label{wcm10}
\Delta_j = 4 \Big( \frac{1}{r^2} C_j - b \Big)^2 \geq 0,
\end{equation}
and so the roots of   (\ref{wcm9}) are real. Since the sum $S_j = y_1 + y_2 = - 2 \omega^2 \leq 0$, the roots are negative
iff $P_j = y_1 y_2 = \omega^4 - \Big( \frac{1}{r^2} C_j - b \Big)^2$ is positive. Taking into account (\ref{wcm1}) and the expression of $b$ we obtain
\begin{equation}\label{wcm11}
P_j = \omega^4 \Big[ 1 - \Big( \frac{12 C_j -(n-1)(n-5)}{6(n-1)} \Big)^2 \Big] 
\end{equation}
where $C_j$ is given by (\ref{s56}).
Thus
$$
P_j = \frac{j(n-j)\Big( j^2-jn+2(n-1) \Big)}{(n-1)^2},
$$
and so $P_j \geq 0$ iff $j^2-jn+2(n-1) \geq 0$ for any $j=0,1,\ldots n-1$. An elementary analysis
leads to 
%
%
%
%
%
%
\begin{theorem}\label{th: sec}
In the $n$-body problem with logarithm interaction the regular $n$-gon relative equilibrium with  is linearly stable iff $n=2,3,4,5,6$.
\end{theorem}

\end{document}